\documentclass[s5,11pt]{article}
\usepackage{amsmath}
\newtheorem{theorem}{Theorem}

\newtheorem {definition}{Definition}

\newtheorem{proposition}{Proposition}[section]
  
\author{ Tord Sj\"odin}
\title{ On Multivariate Hyperbolically Completely Monotone Densities and Their Laplace Transforms }
 \begin{document}  

\maketitle
\begin{abstract} The class HCM consists of all non-negative functions $f$ on $(0,\infty )$ such that $f(uv)\cdot f(u/v)$ is completely monotone with respect to (wrt) $w=v+1/v$, for every fixed $u>0$, and has been extensively studied for a long time. It is closed wrt a number of useful operations, such as pointwise limits and products.  It is also closed wrt products and quotients of independent random variables, some changes of variables and the Laplace transform. We consider its multivariate (bivariate) counterparts MVHCM (BVHCM) and study some of their properties. In particular, we show that MVHCM is closed wrt the Laplace transform and use this to define the class BVGGC-L of bivariate random vectors with Laplace transform in BVHCM. Then BVGGC-L contains Bondesson's class of random vectors in BVGGC in the strong sense and is contained in the corresponding class BVGGC in the weak sense. We also show that BVHCM, in contrast to HCM, is not closed wrt multiplication of independent bivariate random vectors.\end{abstract}
\paragraph{ \it AMS 2010 Subject Classsification:} Primary 60E10
Secondary 62H05
 \section{Introduction.}
 We consider various classes of multivariate probability densities $f(x)$ on the positive cone $R^n_+=(0,\infty )^n$ in $R^n$. In particular, we consider the class of completely monotone functions (CM) and the class of hyperbolically completely monotone functions (HCM), introduced in \cite{B1}, Ch. 5. A function $f$ is completely monotone (CM) if 
$f$ is $C^\infty$ on $R^n_+$ and $(-1)^{|\alpha |}\cdot D^\alpha f(x)\geq 0$, for all  multiindices $\alpha =(\alpha_1,\dots ,\alpha_n)$ of nonnegative integers $\alpha_i$, $1\leq i\leq n$, and $|\alpha |=\alpha_1+\dots +\alpha_n$. The CM-functions on $(0,\infty )$ were characterized by Bernstein \cite{Ber}, p. 56 as Laplace transforms of nonnegative measures. See also Widder \cite{W} Theorem 12 a,b  and Bochner \cite{Boc} Theorem 4.2.1. The class of CM-functions and the related classes of Stieltjes and Bernstein functions are studied in \cite{F} and \cite{SSV}. \\[1em]
The class of hyperbolically completely monotone (HCM) functions on $(0,\infty )$ was introduced in Bondesson \cite{B1} for the study of generalized gamma convolutions (GGC). A function f(x) defined on  $(0,\infty )$ belongs to HCM if $f(uv)\cdot f(u/v)$ is CM in $w=v+1/v$, for all fixed $u>0$. We write $X\sim HCM$ for a random variable $X$, if its density function $f(x)$ is HCM. The class HCM is extensively studied and characterized in \cite{B1} (and there denoted by $\mathcal{C}$).\\ [1em]
The bivariate analogue BVHCM of HCM was introduced in \cite{B0}, p. 193 (and there denoted by $\mathcal C_2$). In contrast to the univariate case, the class BVHCM is much less studied. A comprehensive study of continuous bivariate distributions, without mentioning hyperbolicity, is found in \cite{BL}.\\[1em]
The class MVHCM of n-variate hyperbolically completely monotone functions $f(x)=f(x_1,\dots ,x_n)$ defined in $R^n_+$ and the corresponding class of random vectors is defined in analogy with the bivariate case (see Section 2). We state and prove some of its properties related to marginal and conditional distributions. It was conjectured by Bondesson \cite{B2}, p. 3764, that MVHCM is closed wrt the Laplace transform. We confirm the conjecture in Theorem 1 and use this fact to define the class BVGGC-L of bivariate densitiy functions with Laplace transform in BVHCM. Then BVHCM-L contains Bondesson's class of random vectors $X\sim BVGGC$ in the strong sense (Theorem 2) and is contained in the weak class $BVHCM$. \\[1em]
It is well known that if $X\sim HCM$ and $Y\sim HCM$ are independent random variables, then $X Y\sim HCM$ and $X/Y\sim HCM$, \cite{B1} Theorem 5.1.1. The analogous result for independent random vectors $(X,Y)\sim BVHCM$ and $(X^\prime,Y^\prime )\sim BVHCM$ was posed as an open question in \cite{B2}, p. 3764. Our main result answers this question in the negative (Theorem 3). The reason is that if $f(x,y)$ denotes the density function of the product vector $(XX^\prime,YY^\prime)$, then $f(v_1,v_2)\cdot f(1/v_1,1/v_2)$ depends in a crusial way on the quadratic term $v_1^2/v_2+v_2/v_1^2=w_1\cdot w_3-w_2,$ where as usual $w_1=v_1+1/v_1$, $w_2=v_2+1/v_2$ and $w_3=v_1/v_2+v_2/v_1$. We give an explicit example that exploits this fact.
\section{Preliminaries.} Let $R^n$ denote the $n-$dimensional Euclidean space $R^n$ and consider nonnegative functions $f(x)=f(x_1,\dots ,x_n)$ defined in the positive cone $R^n_+=\{ x\in R^n; x_i>0, \, 1\leq i\leq n\}$ of $R^n$. Integration wrt Lebesgue measure in $R^n$ is written $\int f(x)\,  dx$. Our notation for measures and integrations are standard and for the necessary background in probability theory we rely on \cite{B1} and \cite{F}. The following definition of the class of the n-variate hyperbolically completely monotone functions is in analogy with the bivariate case in \cite{B3}.\begin{definition} A n-variate density $f$ defined on $R^n_+$ is MVHCM  ($MV_nHCM$) if
\begin{equation}f(u_1v_1,u_2v_2,\dots ,u_nv_n)\cdot f(u_1/v_1,u_2/v_2,\dots ,u_n/v_n)\end{equation}
is completely monotone as a function of $w_i=v_i+1/v_i$, $1\leq i\leq n$ and $w_{ij}=v_i/v_j+v_j/v_i$, $1\leq i<j\leq n$, for all positive numbers $u_i$, $1\leq i\leq n$. In the bivariate case $n=2$ we denote $w_{1,2}$ by $w_3$.
\end{definition}
As an example of functions in MVHCM we consider $$f(x)=   \prod _1^n x_i^{\alpha_ i-1}\cdot (1+a_{1,i}x_1 +\cdots +a_{n,i}x_n)^{-\gamma_i},  $$ where $\alpha_i>0,a_i>0$, $1\leq i\leq n$ and $\gamma >0$, properly normalized. Then the product in (1) takes the form
\begin{equation}\prod\limits _{i=1}^n u_i^{2(\alpha_i-1)}\cdot \big(1+\sum\limits _{i=1}^na_iu_i\cdot w_i+ \sum_{i<j}a_ia_ju_iu_j\cdot w_{i,j} )^{-\gamma},   \end{equation}
which clearly is MVHCM, and can be represented as
$$\int \exp \big(- \sum\limits _1^n \lambda _i\cdot w_i-\sum\limits _{i<j}  \lambda_{i,j}w_{i,j}\big)\, d\nu (\lambda_1,\dots , \lambda_m ),
$$
for some nonnegative measure $ \nu$ depending on $a_i, u_i,\alpha_i,\gamma $  and $m={n\choose 2}$. This can be seen directly or by \cite{Boc} Theorem 4.2.1.
\section{Main results.} We begin this section with a number of properties of random vectors $(X,Y)$ with densities $f$ in BVHCM. The first result states that the conditional density $f_{X|Y=y}(x) $, the marginal density $f_X(x)$, and the quotient density $f_{X/Y}(z)$ all are HCM and that BVHVM is closed wrt taking q-th powers, for $|q|\geq 1$.
\begin{proposition} (a) Let $(X,Y)\sim BVHCM$ and denote the density function by $f$. Then $f_{X|Y=y}(x)$, $f_{X}(x)$ and $f_{X/Y}(z)$ are HCM.\\[0.5em]
(b) Let $(X_1,X_2)\sim BVHCM$, then $(X_1^{-1},X_2^{-1})\sim BVHCM$.\\[0.5em]
(c) Let $(X_1,X_2)\sim BVHCM$, then $(X_1^{q},X_2^{q})\sim BVHCM$ for $|q|\geq 1$.
\end{proposition}
{\it Proof.} To prove (a) we first let $F_1(x)=f_{X|Y=y}(x)=f(x,y)$, for some fixed $y>0$. Then
$F_1(uv)\cdot F_1(u/v)=f(uv,y\cdot 1)\cdot f(u/v,y/1)$
is CM in $v+1/v$, since $f$ is BVHCM. Next let $F_2(x)=f_{X}(x)=\int\limits _0^\infty f(x,y)\, dy$, then
$$F_2(uv) F_2(\frac{u}{v})=\int\limits _0^\infty \int\limits _0^\infty f(uv,y)  f(\frac{u}{v},z) \, dydz=2\int\limits _0^\infty \int\limits _0^\infty f(uv,st)  f(\frac{u}{v},\frac{s}{t})\cdot \frac{s}{t}\, dsdt,
$$
by Fubini's Theorem and the hyperbolic change of variables $y=st$, $z=s/t$ with Jacobian $-2s/t$.  For fixed $u$ and $s$ there exists a nonnegative measure $\nu$, depending on $u$ and $s$, such that the last integral equals
$$2\cdot \int\limits _0^\infty s\, ds \int \limits_0^\infty \frac{dt}{t} \int \exp\big(-\lambda _1(v+\frac{1}{v})-\lambda _2(t+\frac{1}{t})-\lambda _3(\frac{v}{t}+\frac{t}{v}) \big)\,     d\nu (\lambda_1,\lambda_2,\lambda_3),$$
 since $f$ is BVHCM. The last two terms in the integrand of the inner most integral, with reversed signs, becomes
$$t\cdot (\lambda_2+\frac{\lambda_3}{v})+\frac{1}{t}\cdot (\lambda_2 +\lambda_3 v)=\rho \cdot \big(\lambda _2^2+\lambda _3^2+\lambda _2\lambda_3\cdot  (v+\frac{1}{v})\big)$$
after a change of variables $t =\rho \cdot (\lambda_2+\lambda_3v)$.  Putting the pieces together, this proves that $F_2$ is CM in the variable $v+1/v$. Finally, let $F_3(z)=f_{X/Y}(z)=\int_0^\infty y\cdot f(zy,y)\, dy$, then $F_3(uv)F_3(u/v)$ equals
$$ \int\limits _0^\infty \int\limits _0^\infty dx dy\, xy \cdot f(uvy,y)\cdot f(\frac{u}{v}x,x)=2\cdot \int\limits _0^\infty \int\limits _0^\infty s^3\, ds\frac{dt}{t}\cdot f(uvst, st)\cdot f(\frac{u}{v}\cdot \frac{s}{t},\frac{s}{t}),$$
by Fubini's Theorem and the same hyperbolic change of variables. The last two factors are CM in the variables $vt+1/vt$, $t+1/t$ and $v+1/v$ and the proof can be completed in analogy with the case $F_2$. \\[0.5em]
For (b), we let $f$ denote the density function of $(X_1,X_2)$ and put $Z=(X_1^{-1},X_2^{-1})$. Then $f_Z(x_1,x_2)=x_1^{-2}x_2^{-2}\cdot f(x_1^{-1},x_2^{-1})$ and $f$ HCM implies that $f_Z$ is HCM.\\[0.5em]
In view of (b), it is enough to prove (c) for $q>1$. Let $f$ be as in (b) and put $Z=(X_1^q,X_2^q)$, then $f_Z(x_1,x_2)=x_1^{1/q-1}x_2^{1/q-1}\cdot f(x_1^{1/q},x_2^{1/q})$. The expression
$$f_Z(u_1v_1,u_2v_2)\cdot f(\frac{u_1}{v_1},\frac{u_2}{v_2})\sim f_Z(u_1^{1/q}v_1^{1/q},u_2^{1/q}v_2^{1/q})\cdot f_Z(u_1^{1/q}v_1^{-1/q},u_2^{1/q}v_2^{-1/q})
$$
can, for fixed $u_1$ and $u_2$, be represented by
$$\int \limits _{R^3_+}\exp \big[ -\lambda_1\cdot  (v_1^{\frac{1}{q}}+v_1^{-\frac{1}{q}})-\lambda_2\cdot  (v_2^{\frac{1}{q}}+ v_2^{-\frac{1}{q}})-\lambda_3\cdot  \big((\frac{v_1}{v_2})^{\frac{1}{q}}+(\frac{v_2}{v_1})^{\frac{1}{q}}\big)\big]\, d\nu ,
$$
for some nonnegative measure $d\nu (\lambda_1,\lambda_2,\lambda_3).$ The terms is the exponent, with reversed sign, are  Bernstein function in the variables $v_1+\frac{1}{v_1}, v_2+\frac{1}{v_2}$ and $\frac{v_1}{v_2}+\frac{v_2}{v_1}$, respectively, by \cite{B0}, p. 60.
\hfill $\triangle$ 
\\[0.5em] {\it Remark 1.} The statements in (a) about $f_{X|Y=y}(x)$ and $f_{X}(x)$ also hold when the random variables $X$ and $Y$ are replaced by random vectors. The details are left to the reader. We show by an example that   $f_{X|Y=y}(x)$, $f_{Y|X=x}(x)$, $f_{X}(x)$ and $f_{Y}(x)$ HCM does not imply that $f_{(X,Y)}$ is BVHCM. \\[0.5em]
{\it Example.} Define $f(x,y)=c\cdot (1+x+y+kxy)^{-\gamma }$, where $\gamma >0$, $c$ is a normalizing constant and $k>1$. Clearly, $f_{X|Y=y}(x)$ and $f_X(x)$ are HCM. Put $H=f(v_1,v_2)\cdot f(1/v_1,1/v_2)$, then
$$H\sim \big(3+k^2+(1+k)(w_1+w_2)+w_3+k(w_1\cdot w_2-w_3)\big)^{-\gamma}$$
and $\frac{\partial H}{\partial w_3}>0$. Thus $f$ is not BVHCM, which completes the example.\hfill$\triangle$\\[1em]
The following multiplication theorem seems to be new and answers a question posed by Bondesson (personal communication).
\begin{proposition}
 Let $(X_1,X_2)\sim BVHCM$  and $Y\sim HCM$, then $Z=(Y\cdot X_1,Y\cdot X_2)\sim BVHCM$.\end{proposition}
Proof. Denote the density functions of $(X_1,X_2)$ and $Y$ by $f(x,y)$ and $g(x)$, respectively. Then
$$f_Z(x,y)=\int\limits _0^\infty f(\frac{x}{t},\frac{y}{t})\cdot g(t)\, \frac{dt}{t^2}$$
and $f_Z(u_1v_1,u_2v_2)\cdot f_Z(u_1/v_1,u_2/v_2)$ becomes
\begin{equation}J=2\cdot \int\limits_0^\infty \int\limits_0^\infty f(\frac{u_1v_1}{xy},\frac{u_2v_2}{xy})\cdot f(\frac{u_1y}{v_1x},\frac{u_2y}{v_2x})\cdot g(xy)\cdot g(\frac{x}{y})\ \, \frac{dx}{x^3}\frac{dy}{y}
\end{equation}
after a hyperbolic change of variables. By our assumptions on $f$ and $g$, this integrand is, for fixed $u_1,u_2$ and $x$, CM in the variables $v_1/y+y/v_1$, $v_2/y+y/v_2$, $v_1/v_2+v_2/v_1$ and $y+1/y$ and can be represented as $\int e^{-E}\, d\nu$, where
$$ E=\lambda_1\cdot (\frac{v_1}{y}+\frac{y}{v_1})+\lambda_2\cdot (\frac{v_2}{y}+\frac{y}{v_2})+\lambda_3\cdot (\frac{v_1}{v_2}+\frac{v_2}{v_1})+\lambda_4\cdot (\frac{ 1}{y}+\frac{y}{ 1})
$$
and $\nu$ is a nonnegative measure. After a change of the order of integration in (3), it is enough to prove that the inner most integral $\int\limits _0^\infty e^{-E}\, \frac{dy}{y}$ is CM. Three of the terms in $E$ contain the variable $y$ and can be written
$$y\cdot (\frac{\lambda_1}{v_1}+\frac{\lambda_2}{v_2}+\lambda_4)+\frac{1}{y}\cdot (\lambda_1v_1+\lambda_2v_2+\lambda_4)=
$$
$$=\rho\cdot (\lambda_1^2+\lambda_2^2+\lambda_4^2+\lambda_1\lambda_2\cdot w_3+\lambda_1\lambda_4\cdot w_1+\lambda_2\lambda_4\cdot w_2)+1/\rho ,$$
after the change of variables $y=\rho \cdot (\lambda_1v_1+\lambda_2v_2+\lambda_4)$, from which it is clear that $J$ is CM and $f_Z$ is BVHCM.
\hfill $\triangle$
\\[1em]Our next result states that the class MVHCM is closed wrt the Laplace transform and answers a question in \cite{B2}, p. 3764. In the univariate case this follows from the multiplication theorem for independent random variables with densities in HCM, see \cite {B1} Theorem 5.1.1.
\begin{theorem}Assume that $f$ is MVHCM, then the Laplace transform of $f$ is also MVHCM. \end{theorem}
Proof. We give the proof for $n=2$ only, the general case being similar. Let $f$ be as in the theorem,  denote the Laplace transform of $f$ by $F$ and define
$H=F(u_1v_1,u_2v_2)\cdot F(\frac{u_1}{v_1},\frac{u_2}{v_2})$. Then
$$H=\int   ds_1\, ds_2\, dt_1 \, dt_2 \,  \exp (-s_1u_1v_1 -s_2u_2v_2-t_1\frac{u_1}{v_1} -t_2\frac{u_2}{v_2}  )\cdot f(s_1,s_2)\cdot f(t_1,t_2 ),$$
where the integration is over $R^4_+$.
Now we make the hyperbolic changes of variables $s_1=xy, t_1=x/y$ and $s_2 =zw, t_2 =z/w$ with Jacobians $-2x/y$ and $-2z/w$ respectively and get
$$H=\int \cdots  \int  dx\, dy\, dz \, dw \,  e^{-u_1v_1xy-u_2v_2zw-\frac{u_1}{v_1}\frac{x}{y}-\frac{u_2}{v_2}\frac{z}{w}}\cdot f(xy,zw)\cdot f(\frac{x}{y},\frac{z}{w} ),$$
which is our starting point.
For fixed $u_1,u_2,x$ and $z$ we get by Bernstein's Theorem, \cite{Boc} Theorem 4.2.1,  that the last two factors can be represented as 
$$ f(xy,zw)\cdot f(\frac{x}{y},\frac{z}{w})=\int e^{-\lambda_1\cdot (y+1/y)-\lambda_2\cdot (w+1/w)-\lambda_3\cdot (y/w+w/y)}\, d\nu (\lambda_1,\lambda_2\lambda_3)
$$
for some nonnegative measure $\nu$ depending on $u_1,u_2,x$ and $z$. Changing the order of integration we find that it is sufficient to prove that the inner most integral
in $H$
$$I= \int \, \frac{dydw}{yw}\,   e^{-u_1x(v_1y+\frac{1}{v_1y})-u_2z(v_2w+\frac{1}{v_2w})
 -\lambda_1\cdot (y+\frac{1}{y})-\lambda_2\cdot (w+\frac{1}{w})-\lambda_3\cdot (\frac{y}{w}+\frac{w}{y})}$$
 is completely monotone as a function of $w_1=v_1+1/v_1$, $w_2=v_2+1/v_2$ and $w_3=v_1/v_2+v_2/v_1$. The exponent in the integrand of $I$ can be written (with the sign changed and new names on the variables)
 $$J=A(v_1y+\frac{1}{v_1y})+B(v_2w+\frac{1}{v_2w})+C (y+\frac{1}{y})+D (w+\frac{1}{w})+E  (\frac{y}{w}+\frac{w}{y})$$
 We are going to perform two changes of variables in $J$ and start with the $y-$variabel. The six terms in $J$ that contain a factor $y$ or $1/y$ can be written 
 $$J_1=  y(Av_1+C+\frac{E}{w})+\frac{1}{y}(\frac{A}{v_1}+C+Ew)$$
 and we make the change of variables $y=\rho\cdot ( A/v_1+C+Ew)$. Then $dt/t=d\rho /\rho$ and  $$J_1= \rho\cdot \big (A^2+C^2+E^2+AC\cdot w_1+CE\cdot (w+\frac{1}{w})+AE\cdot (w v_1+\frac{1}{wv_1})\big)+\frac{1}{\rho} .
$$
We treat the remaining four terms in $K$ together with the terms in  $J_1$ containing a factor $w$ or $1/w$
$$J_2=w\cdot (Bv_2+(D+\rho CE)+\rho AEv_1)+\frac{1}{w}\cdot (\frac{B}{v_2}+(D+\rho CE)+\rho AE\frac{1}{v_1})$$
 in a similar way. The change of variables $w= \delta \cdot (\frac{B}{v_2}+(D+\rho CE)+\rho AE\frac{1}{v_1}), $ $d\delta /\delta =dw/w$, expresses $J_2$  as a linear combination of $w_1$, $w_2$ and $w_3$. Putting the expressions for $J_1$ and $J_2$ together we find that 
 $$J=a+b\cdot w_1+c\cdot w_2+d\cdot w_3,$$
 where a,b,c,d are polynomials in $\lambda_1, \lambda_2, \lambda_3, u_1. u_2, x, z, \rho , 1/\rho , \delta$ and $1/\delta$ with positive coefficients. It follows that $I$, and thereby also $H$, is completely monotone in the variables $w_1$, $w_2$ and $w_3$. This completes the proof.\hfill $\triangle$ 
\noindent  \\[1em]   There is a close connection between HCM and the class of generalized gamma convolutions (GGC) introduced by Thorin \cite{T1}, \cite{T2} in 1977. The class GGC consists of the limit distributions of sums of independent Gamma distributed random variables, see \cite{B1}. Bondesson proved that a density f is GGC if and only its Laplace transform is HCM, \cite{B1} Theorem 5.3.1. Two classes of bivariate GGC were defined in \cite {B1}, pp. 46 -- 47.  A distribution on $R^2_+$ is BVGCC in the strong sense if it is the limit distribution for a sequence of random vectors $Z_n=\big( \sum_{j=1}^n c_{1,j}Y_j,\sum _{j=1}^n c_{2,j}Y_j\big)$, where $Y_j$ are independent, unit scale Gamma variables and $c_{1,j},c_{2,j}\geq 0$. The distribution of a random vector $(X_1,X_2)$, where $X_1,X_2\geq 0$, is BVGGC in the weak sense if $c_1\cdot X_1+c_2\cdot X_2\sim GGC$ for all $c_1,c_2\geq 0$. In view of Theorem 2 and the characterization of GGC mentioned above,  we define BVGGC-L as the class of bivariate densities on $R^2_+$  for which the Laplace transform is BVHCM. Then BVGGC-L is closed wrt addition of random variables and we have the following result.
\begin{theorem}If $X\sim BVGGC$ in the strong sense, then also $X\sim BVGGC$-L.
\end{theorem}
{\it Proof.} Let $Z_n$ be the bivariate random vector in the definition of BVGGC. Then the Laplace transform of $Z_n$ at the point $(s_1,s_2)$ is
$$E\big[ \exp(-\sum\limits _{j=1}^n(s_1c_{1,j}+s_2c_{2,j})Y_j)\big]\sim \prod\limits _{j=1}^n\frac{1}{(1+s_1\cdot c_{1,j}+s_2\cdot  c_{2,j})^{\gamma_j }},
$$
which is BVHCM, since BVHCM is closed wrt multiplication of functions. Taking limits as $n\rightarrow\infty$ completes the proof.\hfill $\triangle$
 \\[1em]It is easy to see that $X\sim BVGGC-L$ implies that $X\sim BVGGC$ in the weak sense, which gives the inclusions
 $$BVGGC_{strong}\subseteq BVGGC-L \subseteq BVGGC_{weak}.$$
  The multiplication theorem for $HCM$ states that if $X\sim HCM$ and $Y\sim HCM$ are independent random variables, then also $X\cdot Y\sim HCM$ and $X/Y\sim HCM$, \cite{B1} Theorem 5.1.1. The corresponding result for bivariate random vectors was posed as an open question in \cite{B2} p. 3764, but so far no proof has appeared. We show by an example that the such a statement is in general false.
  \begin{theorem} There exist independent random vectors $(X,Y)\sim BVHCM $ and $(X^\prime ,Y^\prime ) \sim BVHCM$ such that the product $(X\cdot X^\prime ,Y\cdot Y^\prime )$ does not have a density function in $BVHCM$. In particular, this is the case if $(X,Y) $ and $(X^\prime ,Y^\prime ) $ have densities $f(x,y)=g(x,y)=y^{-2}\cdot e^{-x-k\cdot x/y}$, respectively, for $k>0$ sufficiently large.\end{theorem}
Proof. We set the stage by letting $f$ and $g$ denote the density functions for the two arbitrary independent random vectors $(X,Y)$ and $(X^\prime ,Y^\prime )$, respectively, on $R^2_+$. Since $X$ and $X^\prime$ are independent we get
$$f_{XX^\prime}(x)=\int\limits _0^\infty \frac{ds}{s}\int\limits _0^\infty f(x/s,t)\int\limits _0^\infty g(s,u)\, dtdu$$
and similarily for $f_{YY^\prime}(x)$,
 By the independence assumption, we find the density function for the random vector $(X X^\prime ,Y Y^\prime )$ to be
\begin{equation}F(x,y)=\int\limits _0^\infty \frac{ds}{s}\int\limits _0^\infty \frac{dt}{t}f(x/s,y/t)\cdot g(s,t).
\end{equation}
 Now we put $f(x,y)=g(x,y)=y^{-2}\cdot e^{-x-k\cdot x/y}$ and get
$$F(x,y)=\int\limits _0^\infty\frac{ds}{s}\int\limits _0^\infty\frac{dt}{t}e^{-\frac{x}{s}-k\cdot \frac{xt}{sy}-s-k\cdot \frac{s}{t}}=\int\limits _0^\infty\frac{ds}{s}\int\limits _0^\infty\frac{dt}{t}e^{-s\cdot (1+\frac{k}{t})-\frac{1}{s}\cdot (x+\frac{kxt}{y})}.
$$
From the last formula it can be shown that $F(x,y)$ is finite for all $x,y>0$. Let $v_1,v_2>0$ and define $J=F(v_1,v_2)\cdot F(1/v_1,1/v_2)$, then
$$J=\int\limits _0^\infty\frac{ds}{s}\int\limits _0^\infty\frac{dt}{t}e^{-\frac{v_1}{s}-k\cdot \frac{v_1t}{sv_2}-s-k\cdot \frac{s}{t}}\cdot \int\limits _0^\infty\frac{du}{u}\int\limits _0^\infty\frac{dv}{v}e^{-\frac{1}{v_1u}-k\cdot \frac{v_2v}{v_1u}-u-k\cdot \frac{u}{v}}.$$
We will show that $J$ is not CM in the variables $w_1=v_1+1/v_1, w_2=v_2+1/v_2$ and $w_3=v_1/v_2+v_2/v_1$, provided $k$ is large enough. As a first step we express $J$ as a function of $w_1, w_2$ and $w_3$ by two hyperbolic changes of variables followed by two changes of variables similar to the ones in the proofs of Proposition 3.1 and Theorem 1.
We begin by putting $s=xy,u=x/y$ and $t=zw, v=z/w$, with Jacobians $-2x/y$ and $-2z/w$ respectively, which after some algebra gives $J=4\cdot \int\int\int\int e^{-E}dxdydzdw/xyzw$ and
$$E=\frac{1}{x}\cdot (\frac{v_1}{y}+\frac{y}{v_1})+x\cdot (y+\frac{1}{y})+\frac{kz}{x}\cdot (\frac{v_1w}{v_2y}+\frac{v_2y}{v_1w})+\frac{kx}{z}\cdot (\frac{y}{w}+\frac{w}{y})=$$
$$=y\cdot (\frac{1}{x}v_1+x+k\frac{z}{x}\frac{v_2}{v_1w}+k\frac{x}{z}\frac{1}{w})+\frac{1}{y}\cdot (\frac{1}{x}v_1+x+k\frac{z}{x}\frac{v_1w}{v_2}+k\frac{x}{z}w).
$$
We make a change of variables by putting the second term equal to $1/\rho$ and get
$$E=\rho\cdot \big(\frac{1}{x^2}+x^2+k^2\frac{z^2}{x^2}+k^2\frac{x^2}{z^2}+k\frac{z}{x^2}(\frac{w}{v_2}+\frac{v_2}{w})+\frac{k}{z}(\frac{v_1}{w}+\frac{w}{v_1})+
$$ 
\begin{equation}+kz(\frac{v_2}{v_1w}+\frac{v_1w}{v_2})+k\frac{x^2}{z}(w+\frac{1}{w})+w_1+kw_3\big)+1/\rho .
\end{equation}
 The four terms in (4) that contain the variable $w$ are now treated in the same way. These terms are $\rho$ times
 $$E^\prime =w\cdot \big(\frac{kz}{x^2}\frac{1}{v_2}+\frac{k}{z}\frac{1}{v_1}+kz\frac{v_1}{v_2}+\frac{kx^2}{z}\big)+\frac{1}{w}\cdot \big(\frac{kz}{x^2}v_2+\frac{k}{z}v_1+kz\frac{v_2}{v_1}+\frac{kx^2}{z}\big) .
 $$
 A change of variables by putting the last term equal to $1/\delta$ gives
 $$E^\prime  =\delta\big( \frac{k^2z^2}{x^4}+\frac{k^2x^4}{z^2}+\frac{k^2}{z^2}+k^2z^2+(\frac{k^2z^2}{x^2}+\frac{k^2x^2}{z^2})\cdot w_1+(k^2x^2+\frac{k^2}{x^2})\cdot w_3+
 $$
 $$+k^2\cdot w_2+k^2\cdot (\frac{v_1^2}{v_2}+\frac{v_2}{v_1^2})\big) +1/\delta.$$
 The crusial term here is the nonlinear term $\frac{v_1^2}{v_2}+\frac{v_2}{v_1^2}=w_1w_3-w_2$, which turns out to be responsible for $J$ not being CM. It will be our main concern in the rest of the proof. We also note that $w_2$ is cancelled and that the terms that contain $w_1$ and $w_3$ are of two types. The first type is terms with no $x-$ and $z-$variables, they are $\rho\cdot w_1+\rho k\cdot w_3+\rho \delta k^2\cdot w_1w_3$. The second type of terms are linear combinations of $w_1$ and $w_3$ with coefficients depending on $\rho , \delta , x$ and $z$. If we put all this together we get the following expression for $E$
 $$E=g_1(x,z,\rho ,\delta )+ g_2(w_1, w_3, x, z, \rho , \delta )+\rho\cdot w_1+\rho k\cdot w_3+\rho \delta k^2\cdot w_1w_3,$$
 where $g_1$ and $g_2$ are polynomials in the indicated variables, $1/\rho$ and $1/\delta$. Carrying out the $x-$ and $z-$integration gives the formula
 \begin{equation}J=\int \limits _0^\infty \int \limits _o^\infty  \frac{d\rho}{\rho}\frac{d\delta}{\delta} e^{-1/\rho -\rho /\delta} G (w_1,w_3,\rho ,\delta ) \cdot e^{-\rho\cdot w_1-\rho k\cdot w_3-\rho \delta k^2\cdot w_1w_3}\end{equation}
 on which the remaining part of the proof is based. A careful study of the expression for $E$ shows that $J$ and all its derivatives are finite for $w_1,w_3\geq 0$. We claim that the second order mixed derivative of $J$ is negative for $w_1,w_2$ close to the origin and $k$ sufficiently large, which proves that $J$ is not CM and $F$ is not BVHCM.
  \\[0.5em]
  In the following we denote the integration in (5) by $\int dm$ and suppress the variables $\rho ,\delta$ in $G$. Denoting the derivative of $J$ and $G$ wrt $w_1$ by $J_1$ and $G_1$, respectively, gives
  $$J_1=\int dm\,  e^{-\rho\cdot w_1-\rho k\cdot w_3-\rho \delta k^2\cdot w_1w_3}\cdot \big( G_1(w_1,w_3)-G(w_1,w_3)\cdot (\rho +\rho \delta k^2w_3)\big)
  $$
  and, by letting $w_1\rightarrow 0$,
  $$J_1=\int dm\,  e^{ -\rho k\cdot w_3 }\cdot \big( G_1(0,w_3)-G(0,w_3)\cdot (\rho +\rho \delta k^2w_3)\big).
  $$
  We let $J_{1,3}$ denote the derivative of $J_1$ wrt $w_3$, and similarily for  $G$,  and suppress the $w_1$-variable from the notation. Then
  $$J_{1,3}=\int dm\, e^{-\rho kw_3}\cdot \big[ G_{1,3}(w_3)-G(w_3)\cdot \rho\delta k^2-G_3(w_3)\cdot (\rho +\rho\delta k^2w_3)
  $$
  $$-\rho k\cdot \big(G_1(w_3)-G(w_3)\cdot (\rho +\rho\delta k^2w_3)\big)\big].$$
  Now we let $k\rightarrow\infty$ and $w_3\rightarrow 0$ in such a way that $k^3w_3\rightarrow 0$. Then $J_{1,3}$ is asymptotically
  $$\int dm\, \big( G_{1,3}(0)-G(0)\cdot \rho\delta k^2-G_3(0)-\rho k\cdot G_1(0)+G(0)\cdot \rho^2 k\big) <0,$$
  if $k$ is sufficiently large, since all integrals are finite. This proves our claim and completes the proof of Theorem 3.
   \hfill $\triangle$\\[1em]
{\it Remark 2.} Let $(X,Y)$ and $(X^\prime ,Y^\prime )$ be independent random vectors in BVHCM with density functions $f(x,y) $ and $g(x,y) $, respectively, and define $F(x,y)$ by (4).  Then, in a way similar to the proof of Theorem 3, it can be proved that $J=F(u_1v_1,u_2v_2)\cdot F(u_1/v_1,u_2/v_2)$ equals four times
\begin{equation}    \int\limits _0^\infty \frac{dx}{x}\int\limits _0^\infty \frac{dy}{y}\int\limits _0^\infty \frac{dz}{z}\int\limits _0^\infty \frac{dw}{w}f( \frac{u_1}{x}  \frac{v_1}{y}, \frac{u_2}{z}  \frac{v_2}{w})   
 f( \frac{u_1y}{xv_1} , \frac{u_2w}{zv_2})\cdot   g(xy,zw)  g(\frac{x}{y},\frac{z}{w}).\end{equation}
By assumption, the first two factors and the last two factors in the integrand in (6) can be represented as Laplace transforms of nonnegative measures $\nu_1(\lambda_1,\lambda_2,\lambda_3)$ and $\nu_2(\lambda_4,\lambda_5,\lambda_6)$ respectively. Again, as in the proof of Theorem 3, we can rewrite $J$ as
\begin{equation}J=4\cdot \int\limits_0^\infty\int\limits_0^\infty  e^{-1/\rho-\rho /\delta}\int d\nu_1 \int d\nu_2 \,  e^{-a-b\cdot w_1-c\cdot w_2-d\cdot w_3-e\cdot w_1w_3},
\end{equation} where $a,b,c,d,e$ are polynomials in $\rho ,\delta$ and $ \lambda_1,\dots ,\lambda_6$. The crucial quantity here is the coefficient $e$, which is a product of  $\rho ,\delta$ and four of the $\lambda_i:s.$ If $\nu_1$ and $\nu_2$ have no mass on the set where $e>0$, then $J$ is CM. The converse is probably true, but seems hard to prove in the general case. It follows from (7) that, although $J$ is not in general CM, $J$ is CM in each variable separately.
 \\[1em]{\bf Acknowledgement.} The author thanks Professor Lennart Bondesson for pointing out the problems treated in this paper and illuminating discussions.


\begin{thebibliography}{99}
 \bibitem{BL}Balakrishnan N, Lai C. D, Continuous Bivariate Distributions, Sec. Ed., Springer-Verlag, New York, 2009
 \bibitem{Ber}Bernstein S, {\it Sur les fonctions absolument monotones}, Acta Math. 52, 1929, pp. 1 -- 66
 \bibitem{Boc} Bochner S, Harmonic analysis and the theory of probability, Univ. of Calif. Press, Berkley and Los Angeles, 1960
 \bibitem{B0} Bondesson L, {\it Generalized Gamma Convolutions and Complete Monoticity}, Probability Theory and Related Fields, 85, 1990, pp. 181 -- 194
 \bibitem{B1} Bondesson L,  Generalized Gamma Convolutions and Related Classes of Distributions and Densities, Lecture Notes in Statistics, vol. 76, Springer-Verlag, New York, 1992
 \bibitem{B2}Bondesson L, {\it On univariate and bivariate generalized gamma convolutions}, Journal of Statistical Planning  and Inference, 139, 2009, pp. 3759 -- 3765,  
 \bibitem{B3} Bondesson L, {\it A class of Probability Distributions that are Closed with Respect to Addition as Well as Multiplication of Independent Random Variables}, to appear in Journal of Theoretical Probability
 \bibitem{CTIF}Chatelain F, Tourneret J-Y, Inglada J, Ferrari A, {\it Bivariate Gamma Distributions for Image Registration and Change Detection}, IEEE Trans. on Image Proc. vol, 16, no 7, 2007
 \bibitem{F} Feller W, An introduction to probability theory and its applications, Vol. II, Wiley Series In Probability and Mathematical Statistics, John Wiley and Sons Inc. New York, 1966
 \bibitem{SSV} Schilling R. L, Song R, Vondracek Z,  Bernstein Functions, Studies in Mathematics 37, Walter de Gruyter, Berlin, 2010
 \bibitem{T1} Thorin O, {\it On the infinite divisibility of tha Pareto distribution}, Scand. Actuarial J,. 1977, pp. 31 -- 40
  \bibitem{T2} Thorin O, {\it On the infinite divisibility of tha lognormal distribution}, Scand. Actuarial J,. 1977, pp. 121 -- 148
 \bibitem{W} Widder D. V, {\it Necessary and sufficient conditions for the representation of a function as a Laplace transform}, Trans. Amer. Math. Soc. 33, 1931, pp. 851 -- 892
 \bibitem{W2} Widder D. V, The Laplace transform, Princeton Univ. Press, Princeton, 1946 \end{thebibliography}
 \end{document}